\documentclass{commat}

\title{%
    Stability on Quinquevigintic Functional Equation in Different Spaces
    }

\author{%
    Ramdoss Murali, Sandra Pinelas and Veeramani Vithya
    }

\affiliation{%
    \address{Ramdoss Murali --
    PG and Research Department of Mathematics, Sacred Heart College (Autonomous), Tirupattur
-- 635 601, TamilNadu, India
        }
    \email{%
    shcrmurali@yahoo.co.in
    }
    \address{Sandra Pinelas --
    Departamento de Ciências Exatas e Engenharia, Academia Militar, A. Conde
Castro Guimarães, 2720-113 Amadora, Portugal
        }
    \email{%
    sandra.pinelas@gmail.com
    }
    \address{Veeramani Vithya --
    PG and Research Department of Mathematics, Sacred Heart College (Autonomous), Tirupattur -- 635 601, TamilNadu, India
        }
    }

\abstract{%
    In this work, we determine the general solution of the quinquevigintic
functional equation and also investigate its stability of this equation
in the setting of  matrix normed spaces and the framework of matrix non-Archimedean fuzzy normed spaces
by using the fixed point method.
    }

\keywords{%
     Generalized Hyers-Ulam stability,  fixed point, quinquevigintic functional
equation, matrix normed spaces, matrix non-Archemedian fuzzy normed spaces
    }

\msc{%
    AMS classification 46S40,
46S50, 47L25, 47H10, 54C30, 54E70
    }

\VOLUME{31}
\YEAR{2023}
\NUMBER{1}
\firstpage{221}
\DOI{https://doi.org/10.46298/cm.10341}

\begin{paper} 

\section{Introduction}

In 1940, S. M. Ulam \cite{QDR25} initiated the study of
stability problems for various functional equations. He raised a question
related to the stability of homomorphism. In the following year,
D. H. Hyers \cite{QDR5} was able to give a partial solution to Ulam's
question in the context of Banach spaces. This was the first significant breakthrough and a step towards more studies in this domain of research. In 1978, Th. M. Rassias \cite{QDR12} succeeded in extending the Hyers theorem by considering an unbounded Cauchy difference. He was the first to prove the stability of the linear mapping in Banach spaces. In 1950,  T. Aoki \cite{QDR1} had provided a proof of a special case of the Rassias result when the given function is additive.
In 1994, Gavruta \cite{QDR4} provided a further generalization of Rassias theorem in
which he replaced the bound by a general control function for the existence of a unique linear mapping. Since then, a large number of papers have been published in connection with various generalizations of Ulam problem and Hyers theorem.

In 1991, the fixed point method was used for the first
time by J. A. Baker \cite{bak}, who applied a variant of Banach's fixed point theorem to obtain the Hyers-Ulam stability of a functional equation in a single variable (for more applications of this method \cite{QDR2}, \cite{dia}, \cite{QDR6}, \cite{rad}).

During the last eight decades, the stability problems of various functional equations such as additive \cite{QDR9}, \cite{QDR30}, \cite{add}, quadratic \cite{th}, \cite{quadO4}, cubic \cite{cub}, \cite{cub1}, \cite{QDR17}, quartic \cite{quar2}, \cite{QDR29}, \cite{quar} and mixed types \cite{lee2}, \cite{wan} have been investigated by many mathematicians.

 In the recent years, the stability problems
of higher degree functional equations  (like quintic, sextic, septic, octic upto trevigintic and quattuorvigintic) have been broadly investigated by a number of mathematicians \cite{QDR22}, \cite{mur}, \cite{QDm}, \cite{muru}, \cite{murt}, \cite{QDr}, \cite{QDR34}, \cite{QDR24}, \cite{QDR18}, \cite{QDR28}, \cite{QDR27}. So, we are interested to introducing the new functional equation of the form
\begin{eqnarray}
 &  & \varsigma(u+13v)-25\varsigma(u+12v)+300\varsigma(u+11v)-2300\varsigma(u+10v)+12650\varsigma(u+9v)\nonumber \\
 &  & -53130\varsigma(u+8v)+177100\varsigma(u+7v)-480700\varsigma(u+6v)+1081575\varsigma(u+5v)\nonumber \\
 &  & -2042975\varsigma(u+4v)+3268760\varsigma(u+3v)-4457400\varsigma(u+2v)\nonumber \\
 &  & +5200300\varsigma(u+v)-5200300\varsigma(u)+4457400\varsigma(u-v)-3268760\varsigma(u-2v)\nonumber \\
 &  & +2042975\varsigma(u-3v)-1081575\varsigma(u-4v)+480700\varsigma(u-5v)-177100\varsigma(u-6v)\nonumber \\
 &  & +53130\varsigma(u-7v)-12650\varsigma(u-8v)+2300\varsigma(u-9v)-300\varsigma(u-10v)\nonumber \\
 &  & +25\varsigma(u-11v)-\varsigma(u-12v)=25!\varsigma(v),\label{eq:VT0}
\end{eqnarray}
where the above functional equation (\ref{eq:VT0}) is said to be the quinquevigintic
functional equation if the function $\varsigma(u)=au^{25}$ is its
solution.

The concept of non-Archimedean fuzzy normed spaces and matrix normed spaces has been  introduced by Mirmostafafe et al.
\cite{QDRf14} and Effors et al. \cite{eff} respectively. Quite recently, the new  results on stability
of functional equations in non-Archimedean fuzzy normed spaces and matrix normed spaces studied
in \cite{QDRz12}, \cite{QDR17}, \cite{QDR18} and \cite{QDR7}, \cite{QDR9}, \cite{QDR29} respectively. The first result on the stability of functional equation in the setting of matrix non-Archimedean fuzzy normed spaces has been given in \cite{th}.

In this paper, we determine the general solution of the functional equation (\ref{eq:VT0}) and we also prove the stability of the functional equation (\ref{eq:VT0}) in matrix normed spaces and matrix non-Archimedean fuzzy normed spaces by using fixed point approach.

\section{ The General
Solution of  Functional Equation (\ref{eq:VT0})}

In this part, we provide the general solution of the quinquevigintic functional
equation~(\ref{eq:VT0}). For this, let us consider $\mathcal{D}$
and $\mathcal{E}$ real vector spaces.

\begin{theorem}\label{ge} If $\varsigma:\mathcal{D}\rightarrow\mathcal{E}$
is a mapping satisfying equation~\eqref{eq:VT0} for all $x,y\in\mathcal{D},$
then $\varsigma(2u)=2^{25}\varsigma(u)$ for all $u,v\in\mathcal{D}$.
\end{theorem}

\begin{proof} Letting $u=v=0$ in (\ref{eq:VT0}), one gets $\varsigma(0)=0$.
Substituting $u=0$, $v=u$ and $u=u$, $v=-u$ in (\ref{eq:VT0}) and
adding the two resulting equations, we get $\varsigma(-u)=-\varsigma(u)$.
Hence, $\varsigma$ is an odd mapping.

 Substituting $u=0$, $v=2u$ and $u=12u$, $v=u$ in (\ref{eq:VT0})
and subtracting the two resulting equations, one gets
\begin{align}
 & 25\varsigma(25u)-324\varsigma(24u)+2300\varsigma(23u)-12375\varsigma(22u)+53130\varsigma(21u)-179100\varsigma(20u)\nonumber \\
 & +480700\varsigma(19u)-1071225\varsigma(18u)+2042975\varsigma(17u)-3309240\varsigma(16u)+4457400\varsigma(15u)\nonumber \\
 & -5076330\varsigma(14u)+5200300\varsigma(13u)-4761000\varsigma(12u)+3268760\varsigma(11u)\nonumber \\
 & -1442100\varsigma(10u)+1081575\varsigma(9u)-1442100\varsigma(8u)+177100\varsigma(7u)+1172655\varsigma(6u)\nonumber \\
 & +12650\varsigma(5u)-1190940\varsigma(4u)+300\varsigma(3u)-25!\varsigma(2u)-25!\varsigma(u)=0\label{eq:VT1}
\end{align}
for all $u\in\mathcal{D}.$

Substituting $(u,v)$ by $(12u,u)$ in (\ref{eq:VT0}), and increasing
the resulting equation by $25$, and then subtracting the resulting
equation from (\ref{eq:VT1}), we get
\begin{align}
 & 301\varsigma(24u)-263120\varsigma(21u)+1149150\varsigma(20u)-3946800\varsigma(19u)+10946275\varsigma(18u)\nonumber \\
 & -24996400\varsigma(17u)+47765135\varsigma(16u)-77261600\varsigma(15u)+106358670\varsigma(14u)\nonumber \\
 & -124807200\varsigma(13u)+125246500\varsigma(12u)-108166240\varsigma(11u)+80276900\varsigma(10u)\nonumber \\
 & -49992800\varsigma(9u)-25597275\varsigma(8u)-11840400\varsigma(7u)+5600155\varsigma(6u)-5200\varsigma(23u)\nonumber \\
 & -1315600\varsigma(5u)+874690\varsigma(4u)-57200\varsigma(3u)-25!\varsigma(2u)+25!(26)\varsigma(u)+45125\varsigma(22u)=0\label{eq:VT2}
\end{align}
for all $u\in\mathcal{D}.$

Substituting $(u,v)$ by $(11u,u)$ in (\ref{eq:VT0}), and increasing
the resulting equation by $301$, and then subtracting the resulting
equation from (\ref{eq:VT2}) , we have
\begin{align}
 & 2325\varsigma(23u)-45175\varsigma(22u)+429180\varsigma(21u)-2658500\varsigma(20u)+12045330\varsigma(19u)\nonumber \\
 & -42360825\varsigma(18u)+119694300\varsigma(17u)-277788940\varsigma(16u)+537673875\varsigma(15u)\nonumber \\
 & -1440043800\varsigma(12u)+1457124060\varsigma(11u)-1261400500\varsigma(10u)+933903960\varsigma(9u)\nonumber \\
 & -589338200\varsigma(8u)+313713675\varsigma(7u)-139090545\varsigma(6u)+51991500\varsigma(5u)\nonumber \\
 & -877538090\varsigma(14u)+1216870200\varsigma(13u)-16866820\varsigma(4u)\nonumber \\
 & +3750450\varsigma(3u)-25!\varsigma(2u)+25!(327)\varsigma(u)=0\label{eq:VT3}
\end{align}
for all $u\in\mathcal{D}.$

Substituting $(u,v)$ by $(10u,u)$ in (\ref{eq:VT0}), and increasing
the resulting equation by $2325$, and then  subtracting the resulting
equation from (\ref{eq:VT3}), we have
\begin{align}
 & 12950\varsigma(22u)-268320\varsigma(21u)+2689000\varsigma(20u)-17365920\varsigma(19u)+81166425\varsigma(18u)\nonumber \\
 & -292063200\varsigma(17u)+839838560\varsigma(16u)-1976988000\varsigma(15u)+3872378785\varsigma(14u)\nonumber \\
 & -6382996800\varsigma(13u)+8923411200\varsigma(12u)-1.063357344x10^{10}\varsigma(11u)\nonumber \\
 & +1.0829297x10^{10}\varsigma(10u)-9429551040\varsigma(9u)+7010528800\varsigma(8u)\nonumber \\
 & -4436203200\varsigma(7u)+2375571330\varsigma(6u)-1065636000\varsigma(5u)+394890680\varsigma(4u)\nonumber \\
 & -119776800\varsigma(3u)-25!\varsigma(2u)+25!(2652)\varsigma(u)=0\label{eq:VT4}
\end{align}
for all $u\in\mathcal{D}.$ 

Substituting $(u,v)$ by $(9u,u)$ in (\ref{eq:VT0}),and increasing
the resulting equation by $12950$, and then  subtracting the resulting equation from (\ref{eq:VT4}), we arrive at
\begin{align}
 & 55430\varsigma(21u)-1196000\varsigma(20u)+12419080\varsigma(19u)-82651075\varsigma(18u)+395970300\varsigma(17u)\nonumber \\
 & -1453606440\varsigma(16u)+4248077000\varsigma(15u)-1.013401747x10^{10}\varsigma(14u)\nonumber \\
 & +2.007352945x10^{10}\varsigma(13u)-3.34070308x10^{10}\varsigma(12u)+4.708975656x10^{10}\varsigma(11u)\nonumber \\
 & -5.6514588x10^{10}\varsigma(10u)+5.791433396x10^{10}\varsigma(9u)-5.07128012x10^{10}\varsigma(8u)\nonumber \\
 & +3.78942388x10^{10}\varsigma(7u)-2.4080954x10^{10}\varsigma(6u)+1.294076025x10^{10}\varsigma(5u)\nonumber \\
 & -5830174320\varsigma(4u)+2173655250\varsigma(3u)-25!\varsigma(2u)+25!(15602)\varsigma(u)=0\label{eq:VT5}
\end{align}
for all $u\in\mathcal{D}.$

Substituting $(u,v)$ by $(8u,u)$ in (\ref{eq:VT0}), and increasing
the resulting equation by $55430$, and then  subtracting the resulting equation from (\ref{eq:VT5}), we arrive at
\begin{align}
 & 189750\varsigma(20u)-4209920\varsigma(19u)+44837925\varsigma(18u)-305219200\varsigma(17u)+1491389460\varsigma(16u)\nonumber \\
 & -5568576000\varsigma(15u)+1.651118354x10^{10}\varsigma(14u)-3.98781728x10^{10}\varsigma(13u)\nonumber \\
 & +7.983507345x10^{10}\varsigma(12u)-1.340976102x10^{11}\varsigma(11u)+1.90559094x10^{11}\varsigma(10u)\nonumber \\
 & -2.30338295x10^{11}\varsigma(9u)+2.375398278x10^{11}\varsigma(8u)-2.091794432x10^{11}\varsigma(7u)\nonumber \\
 & -2.091794432x10^{11}\varsigma(7u)+1.571064119x10^{11}\varsigma(6u)-1.00301344x10^{11}\varsigma(5u)\nonumber \\
 & +5.41214725x10^{10}\varsigma(4u)-2.447016x10^{10}\varsigma(3u)-25!\varsigma(2u)+25!(71032)\varsigma(u)=0\label{eq:VT6}
\end{align}
for all $u\in\mathcal{D}.$

Substituting $(u,v)$ by $(7u,u)$ in (\ref{eq:VT0}), and increasing
the resulting equation by $189750$, and then subtracting the resulting equation from (\ref{eq:VT6}), we have
\begin{align}
 & 533830\varsigma(19u)-12087075\varsigma(18u)+131205800\varsigma(17u)-908948040\varsigma(16u)\nonumber \\
 & +4512841500\varsigma(15u)-1.709354147x10^{10}\varsigma(14u)+5.13346522x10^{10}\varsigma(13u)\nonumber \\
 & -1.253937828x10^{10}\varsigma(12u)+2.53556896x10^{10}\varsigma(11u)-4.29688116x10^{11}\varsigma(10u)\nonumber \\
 & +6.15453355x10^{11}\varsigma(9u)-7.492170972x10^{11}\varsigma(8u)+7.775774818x10^{11}\varsigma(7u)\nonumber \\
 & -6.886852381x10^{11}\varsigma(6u)+5.199456763x10^{11}\varsigma(5u)-3.3352829x10^{11}\varsigma(4u)\nonumber \\
 & +1.807017713x10^{11}\varsigma(3u)-25!\varsigma(2u)+25!(260782)\varsigma(u)=0\label{eq:VT7}
\end{align}
for all $u\in\mathcal{D}.$

Substituting $(u,v)$ by $(6u,u)$ in (\ref{eq:VT0}), and increasing
the resulting equation by $533830$, and then subtracting the resulting equation from (\ref{eq:VT7}), we have
\begin{align}
 & 1258675\varsigma(18u)-28943200\varsigma(17u)+318860960\varsigma(16u)+1.630276745x10^{12}\varsigma(8u)\nonumber \\
 & -2240108000\varsigma(15u)+1.126884644x10^{10}\varsigma(14u)-4.32066408x10^{10}\varsigma(13u)\nonumber \\
 & +1.312182982x10^{11}\varsigma(12u)-3.238202862x10^{11}\varsigma(11u)+6.609132283x10^{11}\varsigma(10u)\nonumber \\
 & -1.129508796x10^{12}\varsigma(9u)-1.998498667x10^{12}\varsigma(7u)+2.087390377x10^{12}\varsigma(6u)\nonumber \\
 & -1.85953482x10^{12}\varsigma(5u)+1.411273712x10^{12}\varsigma(4u)\nonumber \\
 & -9.08671764x10^{11}\varsigma(3u)-25!\varsigma(2u)+25!(794612)\varsigma(u)=0\label{eq:VT8}
\end{align}
for all $u\in\mathcal{D}.$

Substituting $(u,v)$ by $(5u,u)$ in (\ref{eq:VT0}), and increasing
the resulting equation by $1258675$, and then  subtracting the resulting equation from (\ref{eq:VT8}), we have
\begin{align}
2523675&\varsigma(17u) - 58741540\varsigma(16u) + 654844500\varsigma(15u) - 4653392315\varsigma(14u) \nonumber \\
&{ }+ 2.366676195x10^{10}\varsigma(13u) - 9.16930443x10^{10}\varsigma(12u) + 2.812247863x10^{11}\varsigma(11u) \nonumber \\
&{ }- 7.004381849x10^{11}\varsigma(10u) + 1.441932762x10^{12}\varsigma(9u) - 2.484029748x10^{12}\varsigma(8u) \nonumber \\
&{ }+ 3.611918019x10^{12}\varsigma(7u) - 4.458065759x10^{12}\varsigma(6u) + 4.68557518x10^{12}\varsigma(5u) \nonumber \\
&{ }- 4.196249281x10^{12}\varsigma(4u) - 3.18971249x10^{12}\varsigma(3u) - 25!\varsigma(2u) + 25!(2053287)\varsigma(u) \nonumber \\
&{ }=0, \label{eq:VT9}
\end{align}
for all $u\in\mathcal{D}.$

Substituting $(u,v)$ by $(4u,u)$ in (\ref{eq:VT0}), and increasing
the resulting equation by $2523675$, and then  subtracting the resulting equation from (\ref{eq:VT9}), we have
\begin{align}
 & 4350335\varsigma(16u)-102258000\varsigma(15u)+1151060185\varsigma(14u)\nonumber \\
 & -8257726800\varsigma(13u)+4.238980845x10^{10}\varsigma(12u)-1.657180562x10^{11}\varsigma(11u)\nonumber \\
 & +5.126923876x10^{11}\varsigma(10u)-1.287611026x10^{12}\varsigma(9u)+2.671772661x10^{12}\varsigma(8u)\nonumber \\
 & -4.63730678\varsigma(7u)+6.790206084x10^{12}\varsigma(6u)-8.43248747x10^{12}\varsigma(5u)\nonumber \\
 & +8.895693333x10^{12}\varsigma(4u)-7.925233602x10^{12}\varsigma(3u)-25!\varsigma(2u)+25!(4576962)\varsigma(u)=0\label{eq:VT10}
\end{align}
for all $u\in\mathcal{D}.$

Substituting $(u,v)$ by $(3u,u)$ in (\ref{eq:VT0}), and increasing
the resulting equation by $4350335$, and then  subtracting the resulting equation from (\ref{eq:VT10}), we have
\begin{align}
 & 6500375\varsigma(15u)-154040315\varsigma(14u)+1748043700\varsigma(13u)-1.26419293x10^{10}\varsigma(12u)\nonumber \\
 & +6.54152423x10^{10}\varsigma(11u)-2.577519409x10^{11}\varsigma(10u)+8.035906583x10^{11}\varsigma(9u)\nonumber \\
 & -2.033332158x10^{12}\varsigma(8u)-4.249013764x10^{12}\varsigma(7u)-7.41998918x10^{12}\varsigma(6u)\nonumber \\
 & +1.090366402x10^{13}\varsigma(5u)-1.349622047x10^{13}\varsigma(4u)\nonumber \\
 & +1.392736917x10^{13}\varsigma(3u)-25!\varsigma(2u)+25!(8927297)\varsigma(u)=0\label{eq:VT11}
\end{align}
for all $u\in\mathcal{D}.$

Substituting $(u,v)$ by $(2u,u)$ in (\ref{eq:VT0}), and increasing
the resulting equation by $6500375$, and then subtracting the resulting equation from (\ref{eq:VT11}), we have
\begin{align}
 & 8469060\varsigma(14u)-202068800\varsigma(13u)+2308933200\varsigma(12u)\nonumber \\
 & -1.681450144x10^{10}\varsigma(11u)+8.760648248x10^{10}\varsigma(10u)-3.474632448x10^{10}\varsigma(9u)\nonumber \\
 & +1.089447992x10^{12}\varsigma(8u)-2.766678464x10^{12}\varsigma(7u)+5.777884692x10^{12}\varsigma(6u)\nonumber \\
 & -9.99913684x10^{12}\varsigma(5u)+1.432733464x10^{13}\varsigma(4u)\nonumber \\
 & -1.675180068x10^{13}\varsigma(3u)-25!\varsigma(2u)+25!(15427672)\varsigma(u)=0\label{eq:VT12}
\end{align}
for all $u\in\mathcal{D}.$

Substituting $(u,v)$ by $(u,u)$ in (\ref{eq:VT0}), and increasing
the resulting equation by $8469060$, and then  subtracting the resulting equation from (\ref{eq:VT12}), we have
\begin{align}
 & 9657700\varsigma(13u)-231784800\varsigma(12u)+2655867500\varsigma(11u)-1.931540001x10^{10}\varsigma(10u)\nonumber \\
 & +9.9957195x10^{10}\varsigma(9u)-3.90943696x10^{11}\varsigma(8u)+1.197265069x10^{12}\varsigma(7u)\nonumber \\
 & -2.93207772x10^{12}\varsigma(6u)+5.803070488x10^{12}\varsigma(5u)-9.28491278x10^{12}\varsigma(4u)\nonumber \\
 & -1.183826379x10^{13}\varsigma(3u)-25!\varsigma(2u)+25!(23896732)\varsigma(u)=0\label{eq:VT13}
\end{align}
for all $u\in\mathcal{D}.$

Substituting $(u,v)$ by $(0,u)$ in (\ref{eq:VT0}), and increasing
the resulting equation by $2496144$, and then  subtracting the resulting equation from (\ref{eq:VT13}), we have $\varsigma(2u)=2^{25}\varsigma(u)$
for all $u\in\mathcal{D}.$ Thus $\varsigma:\mathcal{D}\rightarrow\mathcal{E}$
is a quinquevigintic mapping. \end{proof}

\section{Generalized Hyers-Ulam Stability in Matrix Normed Spaces for Functional Equation (\ref{eq:VT0})}

In this section, we prove the generalized Hyers-Ulam
stability for the functional equation (\ref{eq:VT0}) in matrix normed
spaces by using the fixed point method.

Throughout this section, let us consider $\left(X,\left\Vert .\right\Vert _{n}\right)$
a matrix normed spaces, $\left(Y,\left\Vert .\right\Vert _{n}\right)$
a matrix Banach spaces and let $n$ be a fixed non-negative integer.

For a mapping $\varsigma:X\rightarrow Y$, define $\mathcal{H}\varsigma:X^{2}\rightarrow Y$
and $\mathcal{H}\varsigma_{n}:M_{n}(X^{2})\rightarrow M_{n}(Y)$ by
\begin{align*}
\mathcal{H}\varsigma(c,d) & =\varsigma(c+13d)-25\varsigma(c+12d)+300\varsigma(c+11d)-2300\varsigma(c+10d)\\
 & +12650\varsigma(c+9d)-53130\varsigma(c+8d)+177100\varsigma(c+7d)-480700\varsigma(c+6d)\\
 & +1081575\varsigma(c+5d)-2042975\varsigma(c+4d)+3268760\varsigma(c+3d)-4457400\varsigma(c+2d)\\
 & +5200300\varsigma(c+d)-5200300\varsigma(c)+4457400\varsigma(c-d)-3268760\varsigma(c-2d)\\
 & +2042975\varsigma(c-3d)-1081575\varsigma(c-4d)+480700\varsigma(c-5d)-177100\varsigma(c-6d)\\
 & +53130\varsigma(c-7d)-12650\varsigma(c-8d)+2300\varsigma(c-9d)-300\varsigma(c-10d)\\
 & +25\varsigma(c-11d)-\varsigma(c-12d)-25!\varsigma(d)
\end{align*}
for all $c,d\in X$.
\begin{align*}
\mathcal{H}\varsigma(x_{rs},y_{rs}) & =\varsigma(x_{rs}+13y_{rs})-25\varsigma(x_{rs}+12y_{rs})+300\varsigma(x_{rs}+11y_{rs})\\
 & -2300\varsigma(x_{rs}+10y_{rs})+12650\varsigma(x_{rs}+9y_{rs})-53130\varsigma(x_{rs}+8y_{rs})\\
 & +177100\varsigma(x_{rs}+7y_{rs})-480700\varsigma(x_{rs}+6y_{rs})+1081575\varsigma(x_{rs}+5y_{rs})\\
 & -2042975\varsigma(x_{rs}+4y_{rs})+3268760\varsigma(x_{rs}+3y_{rs})-4457400\varsigma(x_{rs}+2y_{rs})\\
 & +5200300\varsigma(x_{rs}+y_{rs})-5200300\varsigma(x_{rs})+4457400\varsigma(x_{rs}-y_{rs})\\
 & -3268760\varsigma(x_{rs}-2y_{rs})+2042975\varsigma(x_{rs}-3y_{rs})-1081575\varsigma(x_{rs}-4y_{rs})\\
 & +480700\varsigma(x_{rs}-5y_{rs})-177100\varsigma(x_{rs}-6y_{rs})+53130\varsigma(x_{rs}-7y_{rs})\\
 & -12650\varsigma(x_{rs}-8y_{rs})+2300\varsigma(x_{rs}-9y_{rs})-300\varsigma(x_{rs}-10y_{rs})\\
 & +25\varsigma(x_{rs}-11y_{rs})-\varsigma(x_{rs}-12y_{rs})-25!\varsigma(y_{rs})
\end{align*}
for all $x=[x_{rs}],y=[y_{rs}]\in M_{n}(X)$.

\begin{theorem}\label{NO3} Let $q=\pm1$ be fixed and $\sigma:X^{2}\rightarrow[0,\infty)$
be a function such that there exists a $\kappa<1$ with
\begin{equation}
\sigma(c,d)\leq2^{25q}\kappa\sigma(\frac{c}{2^{q}},\frac{d}{2^{q}})\indent\text{\ensuremath{\forall} \ensuremath{c,d\in X}.}\label{QD1}
\end{equation}
Let $\varsigma:X\rightarrow Y$ be a mapping satisfying
\begin{equation}
\left\Vert \mathcal{H}\varsigma_{n}([x_{rs}],[y_{rs}])\right\Vert \leq\sum\limits _{r,s=1}^{n}\sigma(x_{rs},y_{rs})\indent\text{\ensuremath{\forall} \ensuremath{x=[x_{rs}],y=[y_{rs}]\in M_{n}(X)}.}\label{QD2}
\end{equation}
Then there exists a unique quinquevigintic mapping $\mathbb{V}:X\rightarrow Y$
such that
\begin{equation}
\left\Vert \varsigma_{n}([x_{rs}])-\mathbb{V}_{n}([x_{rs}])\right\Vert _{n}\leq\sum\limits _{r,s=1}^{n}\frac{\kappa^{\frac{1-q}{2}}}{2^{25}(1-\kappa)}\sigma^{*}(x_{rs})\indent\text{\ensuremath{\forall} \ensuremath{x=[x_{rs}]\in M_{n}(X)},}\label{QD3}
\end{equation}
where
\begin{eqnarray*}
\sigma^{*}(x_{rs}) & = & \dfrac{1}{25!}[\sigma(0,2x_{rs})+\sigma(13x_{rs},x_{rs})+25\sigma(12x_{rs},x_{rs})+301\sigma(11x_{rs},x_{rs})\\
 &  & +2325\sigma(10x_{rs},x_{rs})+12950\sigma(9x_{rs},x_{rs})+55430\sigma(8x_{rs},x_{rs})\\
 &  & +189750\sigma(7x_{rs},x_{rs})+533830\sigma(6x_{rs},x_{rs})+1258675\sigma(5x_{rs},x_{rs})\\
 &  & +2523675\sigma(4x_{rs},x_{rs})+4350335\sigma(3x_{rs},x_{rs})+6500375\sigma(2x_{rs},x_{rs})\\
 &  & +8469060\sigma(x_{rs},x_{rs})+9657700\sigma(0,x_{rs})].
\end{eqnarray*}
\end{theorem}

\begin{proof} Setting $n=1$ in (\ref{QD2}), we get
\begin{equation}
\left\Vert \mathcal{H}\varsigma(c,d)\right\Vert \leq\sigma(c,d)\label{QD4}.
\end{equation}
Utilizing Theorem \ref{ge}, we get
\begin{align*}
\left\Vert -\varsigma(2c)+2^{25}\varsigma(c)\right\Vert 
\leq \dfrac{1}{25!}&\left[\sigma(0,2c)+\sigma(13c,c)+25\sigma(12c,c)+301\sigma(11c,c)\right.\\
 &{ } +2325\sigma(10c,c)+12950\sigma(9c,c)+55430\sigma(8c,c)\\
 &{ }+189750\sigma(7c,c) +533830\sigma(6c,c) + 1258675\sigma(5c,c) \\
 &{ }+ 2523675\sigma(4c,c) + 4350335\sigma(3c,c) + 6500375\sigma(2c,c) \\
 &{ }\left. + 8469060\sigma(c,c) + 9657700\sigma(0,c) \right].
\end{align*}
Therefore,
\begin{equation}
\left\Vert \varsigma(2c)-2^{25}\varsigma(c)\right\Vert \leq\sigma^{*}(c)\indent\text{\ensuremath{\forall} \ensuremath{c\in X}.}\label{QD39}
\end{equation}
Hence
\begin{equation}
\left\Vert \varsigma(c)-\frac{1}{2^{25q}}\varsigma(2^{q}c)\right\Vert \leq\frac{\kappa^{\left(\frac{1-q}{2}\right)}}{2^{25}}\sigma^{*}(c)\indent\qquad\qquad\text{\ensuremath{\forall} \ensuremath{c\in X.}}\label{QD40}
\end{equation}
Taking $\mathcal{S}=\left\{ f:X\rightarrow Y\right\} $ and
the generalized metric $\rho$ on $\mathcal{S}$ as follows:
\[
\rho(f,g)=\inf\big\{\tau\in\mathbb{R}_{+}:\left\Vert f(c)-g(c)\right\Vert \leq\tau\sigma^{*}(c),\forall c\in X\big\},
\]
it is easy to check that $(\mathcal{S},\rho)$ is a complete generalized
metric (see also \cite{QDR30}). Define the mapping $\mathcal{P}:\mathcal{S}\rightarrow\mathcal{S}$
by
\[
\mathcal{P}f(c)=\dfrac{1}{2^{25q}}f(2^{q}c)\indent\qquad\text{\ensuremath{\forall} \ensuremath{f\in\mathcal{S}} and \ensuremath{c\in X.}}
\]
Let $f,g\in\mathcal{S}$ and $\nu$ an arbitrary
constant with $\rho(f,g)=\nu$. Then

\[
\left\Vert f(c)-g(c)\right\Vert \leq\nu\sigma^{*}(c)
\]
 for all $c\in X.$ Using (\ref{QD1}), we find that

\[
\left\Vert \mathcal{P}f(c)-\mathcal{P}g(c)\right\Vert =\left\Vert \dfrac{1}{2^{25q}}f(2^{q}c)-\dfrac{1}{2^{25q}}g(2^{q}c)\right\Vert \leq\kappa\nu\sigma^{*}(c)
\]
for all $c\in X$. Hence we have
$$
\rho(\mathcal{P}f,\mathcal{P}g)  \leq \kappa\rho(f,g)
$$
for all $f,g\in\mathcal{S}$. By (\ref{QD40}), we have
$$\rho(\varsigma,\mathcal{P}\varsigma)\leq\dfrac{\kappa^{\left(\frac{1-q}{2}\right)}}{2^{25}}.$$
By Theorem 2.2 in \cite{QDR2}, there exists a mapping $\mathbb{V}:X\rightarrow Y$
which satisfies:
\begin{enumerate}
\item $\mathbb{V}$ is a unique fixed point of $\mathcal{P}$, which satisfies $\mathbb{V}(2^{q}c)=2^{{23q}}\mathbb{V}(c)$ $\forall$
$c\in X.$
\item $\rho(\mathcal{P}^{k}\varsigma,\mathbb{V})\rightarrow0$ as $k\rightarrow\infty$.
This implies that $\lim\limits _{k\rightarrow\infty}\frac{1}{2^{25kq}}\varsigma(2^{kq}c)=\mathbb{V}(c),$
$\forall$ $c\in X.$
\item $\rho(\varsigma,\mathbb{V})\leq\dfrac{1}{1-\kappa}\rho(\varsigma,\mathcal{P}\varsigma)$,
which implies the inequality
\begin{equation}
\left\Vert \varsigma(c)-\mathbb{V}(c)\right\Vert \leq\frac{\kappa^{\frac{1-q}{2}}}{2^{25}(1-\kappa)}\sigma^{*}(c)\indent\qquad\qquad\text{\ensuremath{\forall} \ensuremath{c\in X.}}\label{QD26}
\end{equation}

\end{enumerate}
It follows from (\ref{QD1}) and (\ref{QD2}) that
\begin{align*}
\left\Vert \mathcal{H}\mathbb{V}(c,d)\right\Vert
&=\lim\limits _{k\rightarrow\infty}\dfrac{1}{2^{25kq}}\left\Vert \mathcal{H}\varsigma(2^{kq}c,2^{kq}d)\right\Vert \\
&\leq\lim\limits _{k\rightarrow\infty}\dfrac{1}{2^{25kq}}\sigma(2^{kq}c,2^{kq}d) \\
&\leq\lim\limits _{k\rightarrow\infty}\dfrac{2^{kq}\kappa^{q}}{2^{25kq}}\sigma(c,d) \\
&=0,
\end{align*}
 for all $c,d\in X$. Therefore, the mapping $\mathbb{V}:X\rightarrow Y$
is a quinquevigintic mapping. By Lemma 2.1 in \cite{QDR9} and (\ref{QD26}),
we get (\ref{QD3}). Hence $\mathbb{V}:X\rightarrow Y$ is a unique
quinquevigintic mapping satisfying (\ref{QD3}). \end{proof}

\begin{corollary} \label{NO4} Let $q=\pm1$ be fixed and let $l,\omega$
be non-negative real numbers with $l\neq25$. Let $\varsigma:X\rightarrow Y$
be a mapping such that
\begin{equation}
\left\Vert \mathcal{H}\varsigma_{n}([x_{rs}],[y_{rs}])\right\Vert _{n}\leq\sum_{r,s=1}^{n}\omega(\left\Vert x_{rs}\right\Vert ^{l}+\left\Vert y_{rs}\right\Vert ^{l})\indent\text{\ensuremath{\forall} \ensuremath{x=[x_{rs}],y=[y_{rs}]\in M_{n}(X)}.}\label{QD27}
\end{equation}
Then there exists a unique quinquevigintic mapping $\mathbb{V}:X\rightarrow Y$
such that
\[
\left\Vert \varsigma_{n}([x_{rs}])-\mathbb{V}_{n}([x_{rs}])\right\Vert _{n}\leq\sum\limits _{r,s=1}^{n}\dfrac{\omega_{0}}{\left|2^{25}-2^{l}\right|}\left\Vert x_{rs}\right\Vert ^{l}\qquad\forall x=[x_{rs}]\in M_{n}(X),
\]
 where
\begin{align*}
\omega_{0} = \dfrac{\omega}{25!} 
& [34861936 + 6500376(2^{l}) + 4350335(3^{l}) + 2523675(4^{l}) \\
& + 1258675(5^{l}) + 533830(6^{l}) + 189750(7^{l}) + 55430(8^{l}) \\
& + 12950(9^{l}) + 2325(10^{l}) + 301(11^{l}) + 25(12^{l}) + (13^{l})].
\end{align*}
\end{corollary}

\begin{proof} The proof is similar to the proof of Theorem \ref{NO3}
by taking $\sigma(c,d)=\omega(\left\Vert c\right\Vert ^{l}+\left\Vert d\right\Vert ^{l})$
for all $c,d\in X$. Then we can choose $\kappa=2^{q(l-25)}$, and
we obtain the required result. \end{proof}

\section{Quinquevigintic Functional Equation (\ref{eq:VT0}) and its Ulam-Gavruta-Rassias
Stability}
  In this section, we investigate the Ulam-Gavruta-Rassias
stability for the functional equation (\ref{eq:VT0}) in matrix normed
spaces by using the fixed point method.
\begin{theorem} \label{NO5} Let $q=\pm1$ be fixed and let $l,\omega$
be non-negative real numbers such that $l=a+b\neq 25$. Let $\varsigma:X\rightarrow Y$
be a mapping such that
\begin{equation}
\left\Vert \mathcal{H}\varsigma_{n}([x_{rs}],[y_{rs}])\right\Vert _{n}\leq\sum_{r,s=1}^{n}\omega(\left\Vert x_{rs}\right\Vert ^{a}.\left\Vert y_{rs}\right\Vert ^{b})\indent\text{\ensuremath{\forall} \ensuremath{x=[x_{rs}],y=[y_{rs}]\in M_{n}(X)}}.
\end{equation}
Then there exists a unique quinquevigintic mapping $\mathbb{V}:X\rightarrow Y$
such that

\[
\left\Vert \varsigma_{n}([x_{rs}])-\mathbb{V}_{n}([x_{rs}])\right\Vert _{n}\leq\sum\limits _{r,s=1}^{n}\dfrac{\omega_{0}}{\left|2^{25}-2^{l}\right|}\left\Vert x_{rs}\right\Vert ^{l}\qquad\forall x=[x_{rs}]\in M_{n}(X),
\]
 where
\begin{eqnarray*}
\omega_{0} & = & \dfrac{\omega}{25!}[8469060+6500375(2^{a})+4350335(3^{a})+2523675(4^{a})\\
 &  & +1258675(5^{a})+533830(6^{a})+189750(7^{a})+55430(8^{a})\\
 &  & +12950(9^{a})+2325(10^{a})+301(11^{a})+25(12^{a})+(13^{a})].
\end{eqnarray*}
\end{theorem}

\begin{proof} The proof is similar to the proof of Theorem \ref{NO3}.
\end{proof}

\section{Quinquevigintic Functional Equation (\ref{eq:VT0}) and its J. M.
Rassias Stability}
  In this section, we establish the J. M. Rassias
stability for the functional equation (\ref{eq:VT0}) in matrix normed
spaces by using the fixed point method.

\begin{theorem} \label{NO6} Let $q=\pm1$ be fixed and let $l,\omega$
be non-negative real numbers such that $l=a+b\neq 25$. Let $\varsigma:X\rightarrow Y$
be a mapping such that
\begin{equation}
\left\Vert \mathcal{H}\varsigma_{n}([x_{rs}],[y_{rs}])\right\Vert _{n}\leq\sum_{r,s=1}^{n}\omega(\left\Vert x_{rs}\right\Vert ^{a}.\left\Vert y_{ij}\right\Vert ^{b}+\left\Vert x_{rs}\right\Vert ^{a+b}+\left\Vert y_{rs}\right\Vert ^{a+b})
\end{equation}
for all $x=[x_{rs}],y=[y_{rs}]\in M_{n}(X)$. Then there exists
a unique quinquevigintic mapping $\mathbb{V}:X\rightarrow Y$ such
that
\[
\left\Vert \varsigma_{n}([x_{rs}])-\mathbb{V}_{n}([x_{rs}])\right\Vert _{n}\leq\sum\limits _{r,s=1}^{n}\dfrac{\omega_{0}}{\left|2^{25}-2^{l}\right|}\left\Vert x_{rs}\right\Vert ^{l}
\]
for all $x=[x_{rs}]\in M_{n}(X),$ where
\begin{align*}
\omega_{0} = \dfrac{\omega}{25!}
&[ 50492552 + 6500375(2^{a}) + 6500376(2^{l}) + 4350335(3^{a}+3^{l}) \\
&+ 2523675(4^{a}+4^{l}) + 1258675(5^{a}+5^{l}) + 533830(6^{a}+6^{l}) \\
&+ 189750(7^{a}+7^{l}) + 55430(8^{a}+8^{l}) + 12950(9^{a}+9^{l}) \\
&+ 2325(10^{a}+10^{l}) + 301(11^{a}+11^{l}) + 25(12^{a}+12^{l}) + (13^{a}+13^{l})].
\end{align*}
\end{theorem}

\begin{proof} The proof is identical to the proof of Theorem \ref{NO3}.
\end{proof}

\section{Generalized Hyers-Ulam Stability in Matrix\! Non-Archimedean\! Fuzzy
Normed Spaces for (\ref{eq:VT0})}

In this section, we  investigate the generalized Hyers-Ulam
stability for the functional equation (\ref{eq:VT0}) in matrix non-Archimedean
fuzzy normed spaces by using the fixed point method.

Throughout this section, we assume that ${\rm {\mathbb{K}}}$ is a
non-Archimedean field, $X$ is a vector space over ${\rm {\mathbb{K}}}$
and $(Y,{\rm N}_{n})$ is a complete matrix non-Archimedean fuzzy
normed spaces over ${\rm {\mathbb{K}}}$, and $(Z,{\rm N}')$ is (an
Archimedean or a non-Archimedean fuzzy) normed spaces.

\begin{theorem}
Let $q=\pm1$ be fixed and let $\sigma:X\times X\to Z$
be a mapping such that for some $\kappa\neq2^{25}$ with $\left(\frac{\kappa}{2^{25}}\right)^{t}<1$ we have
\begin{equation}
\mathcal{N'}(\sigma(2^{q}c,2^{q}d)\geqq\mathcal{N'}\left(\sigma(c,d),\kappa^{-q}t\right)\label{G4}
\end{equation}
for all $c,d\in X$ and $t>0$, and $\mathop{\lim}\limits _{k\rightarrow\infty}\mathcal{N}(2^{-25kq}\mathcal{H}\varsigma(2^{kq}c,2^{kq}d),t)=1$
for all $c,d\in X$ and $t>0$. Suppose that an odd mapping $\varsigma:X\in Y$
with $\varsigma(0)=0$ satisfies the inequality
\begin{equation}
\mathcal{N}(\mathcal{H}\varsigma_{n}([x_{rs}],[y_{rs}]),t)\geq\sum_{i,j=1}^{n}\mathcal{N}'\left(\sigma(x_{rs},y_{rs}),t\right)\label{G5}
\end{equation}
for all $x=[x_{rs}],y=[y_{rs}]\in M_{n}(X)$, and $t>0$. Then there
exists a unique quinquevigintic mapping $\mathbb{V}:X\to Y$ such
that
\begin{equation}
\mathcal{N}_{n}(\varsigma_{n}([x_{rs}])-\mathcal{V}_{n}([x_{rs}]),t)\geq min\left\{ \Gamma\left(x_{rs},\frac{\left|\kappa-2^{25}\right|t}{n^{2}}\right):r,s=1,2,...,n\right\} \label{G6}
\end{equation}
for all $x=[x_{rs}]\in M_{n}(X)$ and $t>0$, where
\begin{align*}
\Gamma\left(x_{rs},t\right) = \min 
&\left\{ 
\mathcal{N'}\left(\sigma(0,2x_{rs}),\left|25!\right|t\right), \mathcal{N}'\left(\sigma(13x_{rs},x_{rs}),\left|25!\right|t\right)\right.,\\
 & \left.\mathcal{N}'\left(\sigma(12x_{rs},x_{rs}),\frac{\left|25!\right|t}{\left|25\right|}\right), \mathcal{N'}\left(\sigma(10x_{rs},x_{rs}),\frac{\left|25!\right|t}{\left|2325\right|}\right),\right.\\
 &  \left.\mathcal{N}'\left(\sigma(9x_{rs},x_{rs}),\frac{\left|25!\right|t}{\left|12950\right|}\right), \mathcal{N}'\left(\sigma(8x_{rs},x_{rs}),\frac{\left|25!\right|t}{\left|55430\right|}\right),\right.\\
 &  \left.\mathcal{N}'\left(\sigma(7x_{rs},x_{rs}),\frac{\left|25!\right|t}{\left|189750\right|}\right), \mathcal{N}'\left(\sigma(6x_{rs},x_{rs}),\frac{\left|25!\right|t}{\left|533830\right|}\right),\right.\\
 & \left.\mathcal{N}'\left(\sigma(5x_{rs},x_{rs}),\frac{\left|25!\right|t}{\left|1258675\right|}\right), \mathcal{N}'\left(\sigma(4x_{rs},x_{rs}),\frac{\left|25!\right|t}{\left|2523675\right|}\right),\right.\\
 & \left.\mathcal{N}'\left(\sigma(3x_{rs},x_{rs}),\frac{\left|25!\right|t}{\left|4350335\right|}\right), \mathcal{N}'\left(\sigma(2x_{rs},x_{rs}),\frac{\left|25!\right|t}{\left|6500375\right|}\right),\right.\\
 & \mathcal{N}'\left(\sigma(x_{rs},x_{rs}),\frac{\left|25!\right|t}{\left|8469060\right|}\right), \mathcal{N}'\left(\sigma(0,x_{rs}),\frac{\left|25!\right|t}{\left|9657700\right|}\right), \\
 &\left. \mathcal{N'}\left(\sigma(11x_{rs},x_{rs}),\frac{\left|25!\right|t}{\left|301\right|}\right)\right\}.
\end{align*}
\end{theorem}

\begin{proof} For the cases $q=1$ and $q=-1$, we consider $\kappa<2^{25}$
and $\kappa>2^{25}$, respectively. Substituting $n=1$ in (\ref{G5}),
we obtain
\begin{equation}
{\rm {\mathcal{N}}}({\rm {\mathcal{H}}}\varsigma(c,d),t)\ge{\rm {\mathcal{N}}}'\left(\sigma(c,d),t\right)\label{G27}
\end{equation}
for all $c,d\in X$ and $t>0$. Utilizing Theorem 1, we arrive at
\begin{align*}
\mathcal{N} (-25! \varsigma(2c) & + 33554432(25!)\varsigma(c), t) \\
& \ge \min\left\{ \mathcal{N}'\left(\sigma(0,2c),\left|25!\right|t\right),\mathcal{N}'\left(\sigma(13c,c),\left|25!\right|t\right),\right.\\
 &{\qquad} \left.\mathcal{N}'\left(\sigma(10c,c),\frac{\left|25!\right|t}{\left|2325\right|}\right),\mathcal{N}'\left(\sigma(9c,c),\frac{\left|25!\right|t}{\left|12950\right|}\right),\right.\\
 &{\qquad} \left.\mathcal{N}'\left(\sigma(8c,c),\frac{\left|25!\right|t}{\left|55430\right|}\right),\mathcal{N}'\left(\sigma(7c,c),\frac{\left|25!\right|t}{\left|189750\right|}\right),\right.\\
 &{\qquad} \left.\mathcal{N}'\left(\sigma(6c,c),\frac{\left|25!\right|t}{\left|533830\right|}\right),\mathcal{N}'\left(\sigma(5c,c),\frac{\left|25!\right|t}{\left|1258675\right|}\right),\right.\\
 &{\qquad} \left.\mathcal{N}'\left(\sigma(4c,c),\frac{\left|25!\right|t}{\left|2523675\right|}\right),\mathcal{N}'\left(\sigma(3c,c),\frac{\left|25!\right|t}{\left|4350335\right|}\right),\right.\\
 &{\qquad} \left.\mathcal{N}'\left(\sigma(2c,c),\frac{\left|25!\right|t}{\left|6500375\right|}\right),\mathcal{N}'\left(\sigma(c,c),\frac{\left|25!\right|t}{\left|8469060\right|}\right),\right.\\
 &{\qquad} \left.\mathcal{N}'\left(\sigma(0,c),\frac{\left|25!\right|t}{\left|9657700\right|}\right),\mathcal{N}'\left(\sigma(12c,c),\frac{\left|25!\right|t}{\left|25\right|}\right),\right.\\
 &{\qquad} \left.\mathcal{N}'\left(\sigma(11c,c),\frac{\left|25!\right|t}{\left|301\right|}\right)\right\},
\end{align*}
 for all $c\in X$ and $t>0$. Thus
\begin{equation}
{\rm N}\left(\varsigma(c)-\frac{1}{2^{25q}}\varsigma(2^{q}c),\frac{\kappa^{\left(\frac{q-1}{2}\right)}}{\left|2^{25}\right|^{\left(\frac{1+q}{2}\right)}}t\right)\geq\Gamma\left(c,t\right)\label{G32}
\end{equation}
for all $c\in X$ and $t>0$. We consider the set ${\rm {\mathcal{S}_{1}}}=\left\{ f_{1}:X\to Y\right\} $
and introduce the generalized metric $\rho$ on ${\rm {\mathcal{S}_{1}}}$
as follows:
\[
\rho(f_{1},g_{1})=\inf\{\mu\in{\rm R}_{+}:{\rm N}\left(f_{1}(c)-g_{1}(c),t\right)\geq\mu\Gamma\left(c,t\right),\forall c\in X,t>0\}.
\]
It is easy to check that $({\rm {\mathcal{S}_{1}}},\rho)$ is a complete
generalized metric (see Lemma 3.2 in \cite{QDR30}).

 Define the mapping ${\rm {\mathcal{J}}}:{\rm {\mathcal{S}_{1}}}\to{\rm {\mathcal{S}_{1}}}$
by $${\rm {\mathcal{J}}}f_{1}(c)=\dfrac{1}{2^{25q}}f_{1}(2^{q}c)$$
for all $f_{1}\in{\rm {\mathcal{S}_{1}}}$ and $c\in X.$

Let $f_{1},g_{1}\in{\rm {\mathcal{S}_{1}}}$
and $\nu$ be an arbitrary constant with $\rho(f_{1},g_{1})\le\nu$.

Then $${\rm {\mathcal{N}}}\left(f_{1}(c)-g_{1}(c),\nu t\right)\geq\Gamma\left(c,t\right)$$
for all $c\in X$ and $t>0$. 

Therefore, using (\ref{G4}), we get
\[
{\rm {\mathcal{N}}}\left({\rm {\mathcal{J}}}f_{1}(c)-{\rm {\mathcal{J}}}g_{1}(c),\nu t\right)={\rm {\mathcal{N}}}\left(f_{1}(2^{q}c)-g_{1}(2^{q}c),2^{25q}\nu t\right)\geq\Gamma\left(c,\frac{2^{25q}}{\kappa^{q}}t\right)
\]
for all $c\in X$ and $t>0$. Hence by definition $\rho({\rm {\mathcal{J}}}f_{1},{\rm {\mathcal{J}}}g_{1})\le\left(\frac{\kappa}{2^{25}}\right)^{q}\nu$.
This means that ${\rm {\mathcal{J}}}$ is a contractive mapping with
Lipschitz constant $L=\left(\frac{\kappa}{2^{25}}\right)^{l}<1.$
It follows from (\ref{G32}) that
 $$\rho(\varsigma,{\rm {\mathcal{J}}}\varsigma)\le\frac{\kappa^{\left(\frac{q-1}{2}\right)}}{\left|2^{25}\right|^{\left(\frac{1+q}{2}\right)}}.$$
Therefore according to Theorem 2.2 in {[}5{]}, there exists a mapping
$\mathbb{V}:X\to Y$ which satisfies:

\noindent 1. ${\rm {\mathbb{V}}}$ is a unique
fixed point of ${\rm {\mathcal{J}}}$, which satisfies ${\rm {\mathbb{V}}}(2^{q}c)=2^{25q}{\rm {\mathbb{V}}}(c)$
$\forall$ $c\in X.$

\noindent 2. $\rho({\rm {\mathcal{J}}}^{k}\varsigma,{\rm {\mathbb{V}}})\to0$
as $k\to\infty$, which implies that $\mathop{\lim}\limits _{k\to\infty}\frac{1}{2^{25kq}}\varsigma(2^{kq}c)={\rm {\mathbb{V}}}(c)$
$\forall$ $c\in X.$

\noindent 3. $\rho(\varsigma,{\rm {\mathbb{V}}})\le\frac{1}{1-\kappa}\rho(\varsigma,{\rm {\mathcal{J}}}\varsigma)$,
which implies that $\rho(\varsigma,{\rm {\mathbb{V}}})\le\frac{1}{\left|2^{25}-\kappa\right|}.$

So,
\begin{equation}
{\rm {\mathcal{N}}}\left(\varsigma(a)-{\rm {\mathbb{V}}}(c),\frac{1}{\left|2^{25}-\kappa\right|}t\right)\geq\Gamma(c,t)\label{G33}
\end{equation}
for all $c\in X$ and $t>0$. By (\ref{G27}),
\[
{\rm {\mathcal{N}}}({\rm {\mathcal{H}}{\mathbb{V}}}(c,d),t)=\mathop{\lim}\limits _{k\to\infty}{\rm {\mathcal{N}}}(2^{-25kq}{\rm {\mathcal{H}}}\varsigma(2^{kq}c,2^{kq}d),t)\geq\mathop{\lim}\limits _{k\to\infty}{\rm {\mathcal{N}}}'(2^{-25kq}\sigma(2^{kq}c,2^{kq}d),t)=1.
\]
Hence, ${\rm {\mathcal{H}}{\mathbb{V}}}(c,d)$ = 0. Thus, the function
${\rm {\mathbb{V}}}$ satisfies quinquevigintic.

We note that $e_{s}\in M_{1,n}({\rm {\mathbb{R}}})$ means that the
sth  component is $1$ and the others are zero,  $E_{rs}\in M_{n}(X)$
means that (r,s)-component is $1$ and the others are zero, and $E_{rs}\otimes x\in M_{n}(X)$
means that (r,s)-component is $x$ and the others are zero. Since
$\mathcal{N}(E_{kl}\otimes x,t)=\mathcal{N}(x,t)$, we have
\begin{eqnarray*}
\mathcal{N}_{n}([x_{rs}],t) & = & \mathcal{N}_{n}\left(\sum_{r,s=1}^{n}E_{rs}\otimes x_{rs},t\right)\geq min\left\{ \mathcal{N}_{n}(E_{rs}\otimes x_{rs},t_{rs}):\:r,s=1,2,...,n\right\} \\
 & = & \min\left\{ N(x_{rs},t_{rs}):\:r,s=1,2,...,n\right\} .
\end{eqnarray*}
where \quad{}$t=\sum_{r,s=1}^{n}t_{rs}$.

So, $N_{n}([x_{rs}],t)\geq min\left\{ N(x_{rs},\frac{t}{n^{2}}):r,s=1,2,...,n\right\} $.
By (\ref{G33}), we get (\ref{G6}). Thus ${\rm {\mathbb{V}}}:X\to Y$
is a unique quinquevigintic mapping satisfying (\ref{G6}). \end{proof}


\EditInfo{July 05, 2018}{April 30, 2021}{Karl Dilcher}

\end{paper}